\documentclass[12pt]{amsart}

\usepackage{amsmath,amsthm,amsfonts,amssymb}
\usepackage{setspace}
\usepackage[pdftex,pdfborder={0 0 0}]{hyperref}
\usepackage{fullpage}
\usepackage{enumerate}
\usepackage{datetime}
\usepackage[margin=1in]{geometry}
\usepackage{biblatex}
\renewbibmacro{in:}{}
\date{\today}

\newtheorem{innercustomthm}{Theorem}
\newenvironment{customthm}[1]
  {\renewcommand\theinnercustomthm{#1}\innercustomthm}
  {\endinnercustomthm}

\newtheorem{innercustomcorollary}{Corollary}
\newenvironment{customcorollary}[2]
  {\renewcommand\theinnercustomcorollary{#2}\innercustomcorollary}
  {\endinnercustomcorollary}

\newtheorem{innercustomlemma}{Lemma}
\newenvironment{customlemma}[3]
  {\renewcommand\theinnercustomlemma{#2}\innercustomlemma}
  {\endinnercustomlemma}
\title{On Zero-Dimensional Glicci Monomial Ideals}
\author{Benjamin Mudrak}
\date{}

\begin{document}

\maketitle

\begin{abstract} \noindent Consider the polynomial ring $R_n = k[x_1,...,x_n]$, where $k$ is a field. Let $m = (x_1,...,x_n)$ and $I$ be an $m$-primary monomial ideal in $R$. We consider the problem of determining whether such ideals are in the \textbf{G}orenstein \textbf{li}asion \textbf{c}lass of a \textbf{c}omplete \textbf{i}ntersection (\textit{\textbf{glicci}}). We prove that all $m$-primary monomial ideals in $k[x,y,z]$ with at most eight generators are homogeneously glicci. We also construct a large class of $m$-primary monomial ideals in $R_n$ for any $n$ with any number of minimal generators that are homogeneously \textit{glicci} but not in the complete intersection \textbf{li}aison \textbf{c}lass of a \textbf{c}omplete \textbf{i}ntersection (\textit{\textbf{licci}}). All Gorenstein links used are constructed explicitly and every second step links to another $m$-primary monomial ideal.    
\end{abstract}

\section{Introduction}
Let $k$ be a field and $R = k[x_1,...,x_n]$ be the corresponding polynomial ring in $n$ variables. Let $I$ and $J$ be proper ideals. Unless otherwise stated, all ideals in this paper are assumed to be Cohen-Macaulay. We call an ideal generated by a regular sequence a \textit{complete intersection ideal} and say that $I$ and $J$ are \textit{linked} if there exists a complete intersection ideal $C$, such that $C:_RI = J$ and $C:_RJ = I$. Note that if $I$ is unmixed and ${\rm ht}(C) = {\rm ht}(I)$, then $C:_RI = J$ is equivalent to $C:_RJ = I$. When $I$ and $C$ are homogeneous ideals, $I$ and $J$ are said to be \textit{homogeneously} linked. We use the notation $I \sim_C J$ to denote linkage via the complete intersection ideal $C$ and $\mu(I)$ to denote the minimal number of generators of $I$. When the condition is loosened to only require $C$ to be Gorenstein (i.e. $R/C$ a Gorenstein ring), we say that $I$ and $J$ are G-linked. $I$ and $J$ are said to be in the same G-linkage class (resp. linkage class) if there exists a sequence of G-links (resp. links) $I = I_0 \sim_{G_0} I_1 \sim_{G_1} ... \sim_{G_{l-1}} I_l = J$. We also use the term \textit{G-linked} to refer to any two ideals in the same G-linkage class. We call a sequence of two G-links a \textit{double G-link}. Further background on linkage theory can be found in \cite{StructureOfLinkage} and \cite{PeskineSpiro}. 

One major open question in linkage theory concerns identifying which ideals lie in the G-linkage class of a complete intersection ideal. There is a large body of work dedicated to studying the behavior of monomial ideals under linkage. Known classes of glicci monomial ideals include Borel-fixed ideals \cite{BorelFixed}; polarizations of $m$-primary monomial ideals \cite{FK}; and Stanley-Reisner ideals of weakly vertex-decomposable simplicial complexes \cite{GlicciSimplicial}. It is conjectured that every Cohen-Macaulay monomial ideal is homogeneously glicci. Our results provide further evidence supporting this conjecture in the $m$-primary case. We show that this conjecture is true for every $m$-primary monomial ideal in $k[x,y,z]$ that can be generated by eight or fewer elements. We explicitly construct the necessary Gorenstein ideals and every second step links to another $m$-primary monomial ideal.

We also construct a large class of $m$-primary monomial ideals in $k[x_1,...,x_n]$ for any $n$ that are homogeneously glicci despite not being licci. Our explicit construction of every Gorenstein ideal used should allow these ideals to be useful for future experimentation and conjecture-testing.

\section{Preliminaries}
While the glicci case remains open, there is a well-known classification by Huneke and Ulrich of precisely which $m$-primary monomial ideals are licci \cite{HunekeUlrich}. We will make frequent use of this classification, which we restate here for the convenience of the reader:

\begin{customthm}{2.1}{\normalfont  \cite[Theorem 2.6]{HunekeUlrich}}
 Let $R = k[x_1,...,x_n]$ be a polynomial ring over a field $k$. Every $m$-primary monomial ideal $I$ can be written uniquely as $I = (x_1^{a_1},...,x_d^{a_d}) + I^{\#}$, where $I^{\#}$ is generated by monomials that together with $\{x_1^{a_1},...,x_n^{a_n}\}$ generate $I$ minimally. Set $I^{\{0\}} = I$. If $I$ is not a complete intersection then $I^{\#}$ can be written uniquely as $I^{\#} = x^BK$ where $x^B = x_1^{b_1}\cdots x_d^{b_d}$ is a monomial and $K$ a monomial ideal of height at least two. Define $I^{\{1\}} = (x_1^{a_1},...,x_d^{a_d}) + K$ if $I$ is not a complete intersection, and $I^{\{1\}} = S$ otherwise. The following conditions are equivalent:
\begin{enumerate}
    \item $I$ can be linked to a complete intersection by a sequence of links defined by monomial regular sequences.
    \item $I$ can be linked to a complete intersection by a sequence of links defined by homogeneous regular sequences.
    \item $I_{m}$ is licci in $R_{m}$.
    \item $(I^{\{n\}})^{\#}$ has height at most one whenever $I^{\#} \neq R$.
    \item $I^{\{n\}} = R$ for some $n$.
\end{enumerate}
\end{customthm}

\noindent The proof of this result relies on the fact that $I^{\{0\}}$ is homogeneously doubly linked to $I^{\{1\}}$ via the links $(x_1^{a_1},...,x_d^{a_d})$, and $(x_1^{a_1 - b_1},...,x_d^{a_d - b_d})$. Huneke and Ulrich additionally prove that after localizing at $m$, all $m$-primary monomial ideals are glicci [1], but the proof relies heavily on the localization assumption and thus does not carry over to the homogeneous case. In his Ph.D. thesis, Tan Dang was the first to study homogeneous Gorenstein linkage of $m$-primary monomial ideals explicitly and gave a specific subclass of $m$-primary monomial ideals generated by six elements that are homogeneously glicci. We will build off of his work and show it can be greatly generalized.

The analysis of Macaulay Inverse Systems can be a powerful tool for studying $m$-primary Gorenstein ideals. We will briefly review the necessary theory in the graded case that is used in this paper. A similar exposition can be found in \cite{Dang} and further background can be found in \cite{BrunsHerzog}.
\\

\noindent\textbf{Macaulay Inverse Systems:} Let $R = k[x_1,...,x_d]$ and $T = k[x_1^{-1},...,x_d^{-1}]$. $T$ has an explicit $R$-module structure 
defined by the following multiplication operation: Let $x_1^{\alpha_1}\cdots x_d^{\alpha} \in R$ and $x_1^{-\beta_1}\cdots x_d^{-\beta_d} \in T$. Then:

$$x_1^{\alpha_1}\cdots x_d^{\alpha} \cdot x_1^{-\beta_1}\cdots x_d^{-\beta_d} = \begin{cases} 
      x_1^{\alpha_1 - \beta_1}\cdots x_d^{\alpha_d - \beta_d} & \alpha_i - \beta_i \leq 0 \texttt{ }\forall i \\
      0 & \texttt{otherwise}
   \end{cases}
$$
\\

\noindent Through this convention, we associate to any homogeneous $m$-primary $R$-ideal the corresponding annihilator submodule of $T$, ${\rm Ann}(I) = 0 :_{T}I$. This yields a bijective correspondence between finitely generated graded $T$-submodules and homogeneous $m$-primary $R$-ideals. Macaulay Inverse Systems are useful when studying G-linkage due to a well-known result that an $m$-primary homogeneous $R$-ideal $I$ is Gorenstein if and only if ${\rm Ann}(I)$ is cyclic \cite{BrunsHerzog}. The problem of finding $m$-primary homogeneous Gorenstein sub-ideals of $I$ can thus be reduced to finding cyclic graded submodules of $T$ containing ${\rm Ann}(I)$. Macaulay Inverse Systems may also be used to prove the following important double G-linkage construction, which will prove useful in our analysis of homogeneously glicci $m$-primary monomial ideals. We will repeat the proof, so the specific Gorenstein ideals utilized in each double G-link remain clear.

\begin{customthm}{2.2}{\normalfont  \cite[Theorem 5.2.1]{Dang}} Let $R = k[x_1,...,x_n]$ be a polynomial ring over a field $k$, and $I$ be a homogeneous $m$-primary $R$-ideal. If there exists a homogeneous $m$-primary Gorenstein ideal $G$ such that $G \subsetneq I \subset G + (g)$ for some homogeneous element $g \in R$, then $I$ is homogeneously G-linked in two steps to $J = I:_Rg$.
\end{customthm}

\begin{proof}
    Take $F \in T$ such that $0:_TG = RF$ and consider the Gorenstein ideal $G' =0:_RRgF$. It suffices show that $G':_R(I:_R g) = G:_RI$. Now, $G':_R(I:_Rg) = (0:_RRgF):_R(I:_Rg) = 0:_R(I:_Rg)gF= 0:_RIF = (0:_RF):_RI = G:_RI$. Thus $I$ and $I:_Rg$ are homogeneously double G-linked.
\end{proof}

\noindent We apply this theorem to prove our first original result, which we will use throughout the rest of the paper.

\begin{customcorollary}{}{2.3}
    Let $R = k[x_1,...,x_n]$ be a polynomial ring over a field $k$, and consider the ideal $I = (x_1^\alpha,x_2^\beta,x_3^\gamma,x_1^{\alpha_1}x_2^{\beta_1}x_3^{\gamma_1},...,x_1^{\alpha_N}x_2^{\beta_N}x_3^{\gamma_N},x_2^{\beta_{N+1}}x_3^{\gamma_{N+1}})$ where $\alpha \leq \beta_{N+1}+\gamma_{N+1}$. Then $I$ is homogeneously double G-linked to $I = I:_Rx_1^{\min\{\alpha_i\}}$.
\end{customcorollary}

\begin{proof}
    Define $G = (x_2^{\beta},x_3^{\gamma},x_1^\alpha x_2^{\beta - \beta_{N+1}}, x_1^{\alpha}x_3^{\gamma - \gamma_{N+1}},x_1^{\beta_{N+1} + \gamma_{N+1}} - x_2^{\beta_{N+1}}x_3^{\gamma_{N+1}})$. Using an analysis of the corresponding Macaulay Inverse System or the Buchsbaum-Eisenbud structure theorem \cite{BuchsbaumEisenbud}, it is easy to see that $G$ is a Gorenstein ideal. Now, since $\alpha \leq \beta_{N+1}+\gamma_{N+1}$, we have that $G \subsetneq I \subset G + (x_1^{\min{\{\alpha_i\}}})$. Thus by Theorem 2.2, $I$ is homogeneously G-linked to $I :_R x_1^{\min{\{\alpha_i\}}}$.  Notice that this result also extends to $I + (x_{k_1}^{\omega_{k_1}},...,x_{k_l}^{\omega_{k_l}})$ for any $k_1,...,k_l \notin \{1,2,3\}$ and $\omega_{k_1},...,\omega_{k_l} \in \mathbb{N}$. 
\end{proof}

\begin{customthm}{2.4}{\normalfont\cite[Theorem 5.2.2]{Dang}} Let $R = k[x,y,z]$ be a polynomial ring over a field $k$. Let $I$ be an $m$-primary monomial ideal.
\begin{enumerate}
    \item If $\mu(I) \leq 5$, then $I$ is homogeneously licci.
    \item Assume $I = (x^\alpha,y^\beta,z^\gamma,x^{\alpha_1}y^{\beta_1},x^{\alpha_2}z^{\gamma_1},y^{\beta_2}z^{\gamma_2})$ and $\alpha \leq \beta_2 + \gamma_2, \beta \leq \alpha_2 + \beta_1,$ or $\gamma \leq \alpha_1 + \beta_1$. Then $I$ is homogeneously glicci, but not licci.
\end{enumerate}

\end{customthm}

\noindent\textit{Proof of Part 2.} Without loss of generality, assume $\alpha \leq \beta_2 + \gamma_2$. By Corollary 2.3, $I$ is doubly linked to $J = I :_R x^{\min\{\alpha_1,\alpha_2\}}$. Thus $\mu(J) \leq 5$ and so $I$ is homogeneously glicci by Part 1. By Theorem 2.1, this class of ideals is not licci. \qed

\section{Main Results}
To prove our main result, we first show that any nonlicci $m$-primary monomial ideal $I \subset k[x,y,z]$ with $\mu(I) = 6$ is homogeneously G-linked to an ideal of the form described in Theorem 2.4.2. We will then show that any $m$-primary monomial ideal $I \subset k[x,y,z]$ with $\mu(I) = 7$ or $\mu(I) = 8$ can be homogeneously G-linked to an $m$-primary monomial ideal $J$ with $\mu(J) \leq \mu(I) - 1$.

\begin{customlemma}{}{3.1}{} Let $R = k[x,y,z]$ be a polynomial ring over a field $k$. Let $I$ be an $m$-primary monomial $R$-ideal with $\mu(I) \leq 6$. Then $I$ is homogeneously glicci via a sequence of G-links with monomial ideals in each second step.
\end{customlemma}

\begin{proof}

Let $I \subset k[x,y,z]$ be an $m$-primary monomial ideal with $\mu(I) \leq 6$. If $\mu(I) \leq 5$, then $I$ homogeneously glicci by Theorem 2.4.2. Clearly, if $I$ is homogeneously licci, then $I$ is homogeneously glicci. If $I$ is not homogeneously licci, then  $I$ is homogeneously linked to an ideal of the form $(x^\alpha,y^\beta,z^\gamma,x^{\alpha_1}y^{\beta_1},x^{\alpha_2}z^{\gamma_1},y^{\beta_2}z^{\gamma_2})$ by Theorem 2.1. Thus we may assume $I = (x^\alpha,y^\beta,z^\gamma,x^{\alpha_1}y^{\beta_1},x^{\alpha_2}z^{\gamma_1},y^{\beta_2}z^{\gamma_2})$ for some positive integers $\alpha \geq \max\{\alpha_1,\alpha_2\}, \beta \geq \max\{\beta_1,\beta_2\}, $ and $\gamma \geq \max\{\gamma_1,\gamma_2\}$. If $\alpha \leq \beta_2 + \gamma_2, \beta \leq \alpha_2 + \gamma_2,$ or $\gamma \leq \alpha_1 + \beta_1$, then $I$ is homogeneously glicci by Theorem 2.4. Thus we may assume $\alpha > \beta_2 + \gamma_2, \beta > \alpha_2 + \gamma_2,$ and $\gamma > \alpha_1 + \beta_1$.  Take $n \in \mathbb{N}$ with $n > \alpha+\beta+\gamma$. Define $I' = (x^{n},y^{n},z^{n}):_RI = (x^{n},y^{n},z^{n},x^{n - \alpha_1}y^{n - \beta}z^{n - \min(\gamma_i)},x^{n - \text{min}(\alpha_i)}y^{n - \beta_2}z^{n - \gamma},x^{n - \alpha_1}y^{n - \beta_2}z^{n - \gamma_1},x^{n - \text{min}(\alpha_i)}y^{n - \beta}z^{n-\gamma_2},\\x^{n-\alpha}y^{n-\beta_1}z^{n - \text{min}(\gamma_i)},x^{n - \alpha_2}y^{n - \text{min}(\beta_i)}z^{n - \gamma},x^{n-\alpha}y^{n-\text{min}(\beta_i)}z^{n-\gamma_1},x^{n-\alpha_2}y^{n-\beta_1}z^{n-\gamma_2})$. Consider the Gorenstein ideal $G = (x^{n},y^{n},z^{n})$. Then we have $G \subsetneq I' \subset G + (x^{n-\alpha}y^{n-\beta}z^{n-\gamma})$. Thus, by Theorem 2.2, $I'$ is homogeneously double G-linked to $I_1 = I' :_R (x^{n-\alpha}y^{n-\beta}z^{n-\gamma}) = (x^\alpha,y^\beta,z^\gamma,x^{\alpha - \alpha_1}z^{\gamma-\text{min}(\gamma_i)},x^{\alpha - \text{min}(\alpha_i)}y^{\beta-\beta_2},x^{\alpha-\alpha_1}y^{\beta-\beta_2}z^{\gamma-\gamma_1},x^{\alpha-\text{min}(\alpha_i)}z^{\gamma-\gamma_2},y^{\beta-\beta_1}z^{\gamma-\text{min}(\gamma_i)},\\x^{\alpha-\alpha_2}y^{\beta-\text{min}(\beta_i)},y^{\beta-\text{min}(\beta_i)}z^{\gamma-\gamma_1},x^{\alpha-\alpha_2}y^{\beta-\beta_1}z^{\gamma-\gamma_2})$. 
\\

\noindent Write $I_1 = (x^{\alpha
},y^\beta,z^\gamma) + I_1^\# = (x^\alpha,y^\beta,z^\gamma) + x^{\alpha - \max\{\alpha_i\}}K_x + K_{yz}$, where $K_{yz} \subset (yz)$, $K_x \subset (x)$, and $K_{yz} :_R (x) = K_{yz}$. Similarly, define $K_{xz}$ and $K_{xy}$ to be the sub-ideals of $I_1^\#$ such that $K_{xy} \subset (xy)$, $K_{xz} \subset (xz)$, $K_{xy}:_R(z) = K_{xy}$, and $K_{xz}:_R(y) = K_{xz}$. Clearly, $K_{yz} = (y^{\beta-\text{min}(\beta_i)}z^{\gamma - \gamma_1},y^{\beta - \beta_1}z^{\gamma-\text{min}(\gamma_i)}), K_{xy} = (x^{\alpha-\alpha_2}y^{\beta-\text{min}(\beta_i)},x^{\alpha - \text{min}(\alpha_i)}y^{\beta - \beta_2}), \text{ and }K_{xz} = (x^{\alpha -\text{min}(\alpha_i)}z^{\gamma-\gamma_2},x^{\alpha - \alpha_1}z^{\gamma -\text{min}(\gamma_i)})$. Define $T(I) = (i,j,k)$ where $\alpha_i = \text{min}(\alpha_i)$, $\beta_j = \text{min}(\beta_j)$, and $\gamma_k = \text{min}(\gamma_k)$. Notice that if $T(I) = (1,2,1)$ or $T(I) = (2,1,2)$ then $K_{xy},K_{xz},$ and $K_{yz}$ are cyclic. $T(I)$ turns out to be useful for determining which homogeneous G-links are needed to link $I_1$ to a complete intersection. Let $\mathcal{T}$ be the set of possible values of $T(J)$ for some $m$-primary monomial ideal $J$. Notice that if $T(I), (a,b,c) \in \mathcal{T}\setminus\{(2,1,2),(1,2,1)\}$, then we may re-label the variables such that $T(I) = (a,b,c)$. Similarly, if $T(I) = (1,2,1)$, then we may re-label the variables to get $T(I) = (2,1,2)$ and vice versa. Thus we need only consider two cases based on the value of $T(I)$.
\\

\noindent \textbf{Case 1:} Say $T(I) = (1,2,1)$. Thus $I_1 = (x^\alpha,y^\beta,z^\gamma,x^{\alpha - \alpha_2}y^{\beta - \beta_2},y^{\beta - \beta_1}z^{\gamma-\gamma_1},x^{\alpha-\alpha_1}z^{\gamma-\gamma_2},\\x^{\alpha-\alpha_2}y^{\beta-\beta_1}z^{\gamma-\gamma_2})$. Notice that either $\alpha \leq \beta - \beta_1 + \gamma - \gamma_1$, $\beta \leq \alpha - \alpha_1 + \gamma - \gamma_2$, or $\gamma \leq \alpha - \alpha_2 + \beta - \beta_2$. To see this, assume for sake of contradiction that $\alpha > \beta - \beta_1 + \gamma - \gamma_1$, $\beta > \alpha - \alpha_1 + \gamma - \gamma_2$, and $\gamma > \alpha - \alpha_2 + \beta - \beta_2$. Thus $\alpha + \beta + \gamma > 2\alpha + 2\beta + 2\gamma - \alpha_1 - \alpha_2 - \beta_1 - \beta_2 - \gamma_1 - \gamma_2$. From our earlier assumption, we also know that $\alpha > \beta_2 + \gamma_2$, $\beta > \alpha_2 + \gamma_1$, and $\gamma > \alpha_1 + \beta_1$. Thus  $\alpha + \beta + \gamma > \alpha_1 + \alpha_2 + \beta_1 + \beta_2 + \gamma_1 + \gamma_2$. Thus $2\alpha + 2\beta + 2\gamma - \alpha_1 - \alpha_2 - \beta_1 - \beta_2 - \gamma_1 - \gamma_2 > 2\alpha + 2\beta + 2\gamma - \alpha - \beta - \gamma = \alpha + \beta + \gamma$. Thus $\alpha + \beta + \gamma > 2\alpha + 2\beta + 2\gamma - \alpha_1 - \alpha_2 - \beta_1 - \beta_2 - \gamma_1 - \gamma_2 > \alpha + \beta + \gamma$, a contradiction. Without loss of generality, we may assume that $\alpha \leq \beta - \beta_1 + \gamma - \gamma_1$. Then $I_1$ is double G-linked to $I_2 = I_1 :_R (x^{\alpha - \alpha_2})$
$ = (x^\alpha,y^\beta,z^\gamma,x^{\alpha - \alpha_2}y^{\beta - \beta_2},y^{\beta - \beta_1}z^{\gamma-\gamma_1},x^{\alpha-\alpha_1}z^{\gamma-\gamma_2},x^{\alpha-\alpha_2}y^{\beta-\beta_1}z^{\gamma-\gamma_2}):_R (x^{\alpha - \alpha_2})$
$= (x^{\alpha_2},y^{\beta-\beta_2},z^\gamma,x^{\alpha_2-\alpha_1}z^{\gamma-\gamma_2},y^{\beta-\beta_1}z^{\gamma-\gamma_2})$. Thus $\mu(I_2) \leq 5$ and $I$ is homogeneously glicci by Theorem 2.4.1.
\\

\noindent \textbf{Case 2:} Say $T(I) = (1,1,1)$. As before, $\alpha > \beta_2 + \gamma_2, \beta > \alpha_1 + \beta_1, \text{ and }\gamma > \alpha_1 + \beta_1$. Notice that for $n \geq \text{max}(\alpha,\beta+\alpha_1-\gamma+\gamma_2)$, 
$I_1 = (x^n,y^\beta,z^\gamma) :_R I = (x^n,y^\beta,z^\gamma,x^{n-\alpha_1}y^{\beta-\beta_2},x^{n-\alpha_1}z^{\gamma-\gamma_2},\\x^{n-\alpha}y^{\beta-\beta_1}z^{\gamma-\gamma_1},x^{n-\alpha_2}y^{\beta-\beta_1})$. Since $\beta \leq n - \alpha_1 + \gamma - \gamma_2$, $I_1$ is double G-linked to $I_2 = I_1 :_R (y^{\beta-\beta_2}) = 
(x^n,y^{\beta_2},z^\gamma,x^{n-\alpha_1},x^{n-\alpha_1}z^{\gamma-\gamma_2},x^{n-\alpha}y^{\beta_2-\beta_1}z^{\gamma-\gamma_1},x^{n-\alpha_2}y^{\beta_2-\beta_1})
= (x^{n-\alpha_1},y^{\beta_2},z^\gamma,\\x^{n-\alpha}y^{\beta_2-\beta_1}z^{\gamma-\gamma_1},x^{n-\alpha_2}y^{\beta_2-\beta_1})$. Thus $\mu(I_2) \leq 5$ and $I$ is homogeneously glicci.  We note that every second constructed ideal was $m$-primary and monomial.
\end{proof}

\begin{customlemma}{}{3.2}{}
Let $R = k[x,y,z]$ be a polynomial ring over a field $k$. Let $I$ be an $m$-primary monomial $R$-ideal with $\mu(I) \leq 7$. Then $I$ is homogeneously glicci via a sequence of G-links with monomial ideals in each second step.
\end{customlemma}

\begin{proof}
Let $I \subset k[x,y,z]$ be an $m$-primary monomial ideal with $\mu(I) \leq 7$. If $\mu(I) \leq 6$ then $I$ is homogeneously glicci by Lemma 3.1. Similarly, if $I$ is homogeneously licci, then $I$ is homogeneously glicci. If $I$ is not homogeneously licci, then $I$ is homogeneously linked to an ideal of the form $(x^\alpha,y^\beta,z^\gamma,x^{\alpha_1}y^{\beta_1},x^{\alpha_2}z^{\gamma_1},y^{\beta_2}z^{\gamma_2},x^{\alpha_3}y^{\beta_3}z^{\gamma_3})$ by Theorem 2.1. Thus without loss of generality we may assume $I = (x^\alpha,y^\beta,z^\gamma,x^{\alpha_1}y^{\beta_1},x^{\alpha_2}z^{\gamma_1},y^{\beta_2}z^{\gamma_2},x^{\alpha_3}y^{\beta_3}z^{\gamma_3})$. Take $n \in \mathbb{N}$ such that $\gamma < n - \alpha + \alpha_1 + n - \beta + \beta_1$. Then I is homogeneously double G-linked to $I_2 = ((x^\alpha,y^{n},z^\gamma):_RI):_R(y^{n-\beta}) = ((x^\alpha,y^{n},z^\gamma) :_R(x^\alpha,y^\beta,z^\gamma,x^{\alpha_1}y^{\beta_1},x^{\alpha_2}z^{\gamma_1},y^{\beta_2}z^{\gamma_2},x^{\alpha_3}y^{\beta_3}z^{\gamma_3})):_R(y^{n-\beta})= (x^\alpha,y^\beta,z^\gamma,x^{\alpha-\alpha_3}y^{\beta-\beta_2},x^{\alpha-\alpha_2}y^{\beta-\beta_3},y^{\beta-\beta_3}z^{\gamma-\gamma_1},y^{\beta-\beta_1}z^{\gamma-\gamma_3},x^{\alpha-\alpha_1}z^{\gamma-\gamma_3},x^{\alpha-\alpha_3}z^{\gamma-\gamma_2})$ which in turn is homogeneously double G-linked to $I_3 = ((x^{n},y^{n},z^{n}):_RI_2):_R(z^{n - \gamma}) = $ $(x^{n},y^{n},z^{\gamma},x^{n - \alpha + \alpha_3}y^{n - \beta + \beta_3}z^{\gamma_3},x^{n - \alpha + \alpha_1}y^{n - \beta + \beta_1},x^{n - \alpha + \alpha_2}y^{n - \beta}z^{\gamma_1},x^{n - \alpha}y^{n - \beta + \beta_2}z^{\gamma_2})$ by Corollary 2.3.
$I_3$ is then homogeneously double G-linked to $I_4 = (I_3 :_R (z^{\min\{\gamma_i\}})) = (x^{n},y^{n},z^{\gamma-\min\{\gamma_i\}},\\x^{n - \alpha + \alpha_3}y^{n - \beta + \beta_3}z^{\gamma_3 - \min\{\gamma_i\}}, x^{n - \alpha + \alpha_1}y^{n - \beta + \beta_1}, x^{n - \alpha + \alpha_2}y^{n - \beta}z^{\gamma_1 - \min\{\gamma_i\}},x^{n - \alpha}y^{n - \beta + \beta_2}z^{\gamma_2 - \min\{\gamma_i\}})$.\\ We will now break into cases based off of $\min\{\gamma_i\}$.
\\

\noindent\textbf{Case 1:} If $\min\{\gamma_i\} = \gamma_3$, then $I_4 = (x^{n},y^{n},z^{\gamma-\gamma_1},x^{n - \alpha + \alpha_3}y^{n - \beta + \beta_3}z^{\gamma_3 - \gamma_1}, x^{n - \alpha + \alpha_1}y^{n - \beta + \beta_1},\\ x^{n - \alpha + \alpha_2}y^{n - \beta},x^{n - \alpha}y^{n - \beta + \beta_2}z^{\gamma_2 - \gamma_1})$. $I_4$ is then homogeneously double G-linked to $I_5 = I_4 :_R (x^{n - \alpha}) = (x^{\alpha},y^{n},z^{\gamma-\gamma_1},x^{\alpha_3}y^{n - \beta + \beta_3}z^{\gamma_3 - \gamma_1}, x^{\alpha_1}y^{n - \beta + \beta_1}, x^{\alpha_2}y^{n - \beta},y^{n - \beta + \beta_2}z^{\gamma_2 - \gamma_1})$ which is homogeneously double G-linked to $I_6 = I_5 :_R (y^{n - \beta})= (x^{\alpha},y^{\beta},z^{\gamma-\gamma_1},x^{\alpha_3}y^{\beta_3}z^{\gamma_3 - \gamma_1}, x^{\alpha_1}y^{\beta_1}, x^{\alpha_2},\\y^{\beta_2}z^{\gamma_2 - \gamma_1}) = (x^{\alpha_2},y^{\beta},z^{\gamma-\gamma_1},x^{\alpha_3}y^{\beta_3}z^{\gamma_3 - \gamma_1}, x^{\alpha_1}y^{\beta_1}, y^{\beta_2}z^{\gamma_2 - \gamma_1})$. Notice that $\mu(I_6) \leq 6$, thus $I_6$ is homogeneously glicci and thus $I$ is homogeneously glicci. 
\\

\noindent\textbf{Case 2:} If $\min\{\gamma_i\} = \gamma_2$, then $I_4 = (x^{n},y^{n},z^{\gamma-\gamma_2},x^{n - \alpha + \alpha_3}y^{n - \beta + \beta_3}z^{\gamma_3 - \gamma_2}, x^{n - \alpha + \alpha_1}y^{n - \beta + \beta_1}, \\x^{n - \alpha + \alpha_2}y^{n - \beta}z^{\gamma_1 - \gamma_2},x^{n - \alpha}y^{n - \beta + \beta_2})$. Which by Corollary 2.3 is homogeneously double G-linked to $I_5 = I_4 :_R (x^{n - \alpha}) = (x^{\alpha},y^{n},z^{\gamma-\gamma_2},x^{\alpha_3}y^{n - \beta + \beta_3}z^{\gamma_3 - \gamma_2}, x^{\alpha_1}y^{n - \beta + \beta_1}, x^{\alpha_2}y^{n - \beta}z^{\gamma_1 - \gamma_2},y^{n - \beta + \beta_2})$ which in turn is homogeneously double G-linked to $I_6 = I_5 :_R (y^{n - \beta})= (x^{\alpha},y^{\beta},z^{\gamma-\gamma_2},\\x^{\alpha_3}y^{\beta_3}z^{\gamma_3 - \gamma_2}, x^{\alpha_1}y^{\beta_1}, x^{\alpha_2}z^{\gamma_1 - \gamma_2},y^{\beta_2}) = (x^{\alpha},y^{\beta_2},z^{\gamma-\gamma_2},x^{\alpha_3}y^{\beta_3}z^{\gamma_3 - \gamma_2}, x^{\alpha_1}y^{\beta_1}, x^{\alpha_2}z^{\gamma_1 - \gamma_2})$. Notice that $\mu(I_6) \leq 6$ and thus $I_6$ is homogeneously glicci and thus $I$ is homogeneously glicci.
\\

\noindent\textbf{Case 3:} If $\min\{\gamma_i\} = \gamma_3$, then $I_4 = (x^{n},y^{n},z^{\gamma-\gamma_3},x^{n - \alpha + \alpha_3}y^{n - \beta + \beta_3}, x^{n - \alpha + \alpha_1}y^{n - \beta + \beta_1},\\ x^{n - \alpha + \alpha_2}y^{n - \beta}z^{\gamma_1 - \gamma_3},x^{n - \alpha}y^{n - \beta + \beta_2}z^{\gamma_2 - \gamma_3})$. $I_4$ is then homogeneously double G-linked to $I_5 = I_4 :_R (x^{n - \alpha}) = (x^\alpha,y^{n},z^{\gamma-\gamma_3},x^{\alpha_3}y^{n - \beta + \beta_3},x^{\alpha_1}y^{n - \beta + \beta_1},x^{\alpha_2}y^{n - \beta}z^{\gamma_1 - \gamma_3},y^{n - \beta + \beta_2}z^{\gamma_2 - \gamma_3})$ which in turn is homogeneously double G-linked to $I_6 = I_5 :_R (x^{\min(\alpha_i)})= (x^{\alpha-\min(\alpha_i)},y^{n},z^{\gamma-\gamma_3},\\x^{\alpha_3 - \min(\alpha_i)}y^{n - \beta + \beta_3},x^{\alpha_1 - \min(\alpha_i)}y^{n - \beta + \beta_1},x^{\alpha_2 - \min(\alpha_i)}y^{n - \beta}z^{\gamma_1 - \gamma_3},y^{n - \beta + \beta_2}z^{\gamma_2 - \gamma_3})$. We finally divide again into sub-cases based on the value of $\min\{\alpha_i\}$.
\\

\noindent\textbf{Subcase a):} If $\min\{\alpha_i\} = \alpha_1$, then $I_6 = (x^{\alpha-\alpha_1},y^{n},z^{\gamma-\gamma_3},x^{\alpha_3 - \alpha_1}y^{n - \beta + \beta_3},y^{n - \beta + \beta_1},\\x^{\alpha_2 - \alpha_1}y^{n - \beta}z^{\gamma_1 - \gamma_3},y^{n - \beta + \beta_2}z^{\gamma_2 - \gamma_3}) = (x^{\alpha-\alpha_1},y^{n},z^{\gamma-\gamma_3},x^{\alpha_3 - \alpha_1}y^{n - \beta + \beta_3},x^{\alpha_2 - \alpha_1}y^{n - \beta}z^{\gamma_1 - \gamma_3},\\y^{n - \beta + \beta_2}z^{\gamma_2 - \gamma_3})$. Notice that $\mu(I_6) \leq 6$, thus $I_6$ is homogeneously glicci and thus $I$ is homogeneously glicci.
\\

\noindent\textbf{Subcase b):} If $\min\{\alpha_i\} = \alpha_2$, then $I_6 = (x^{\alpha-\alpha_2},y^{n},z^{\gamma-\gamma_3},x^{\alpha_3 - \alpha_2}y^{n - \beta + \beta_3},x^{\alpha_1 - \alpha_2}y^{n - \beta + \beta_1},\\y^{n - \beta}z^{\gamma_1 - \gamma_3},y^{n - \beta + \beta_2}z^{\gamma_2 - \gamma_3})$ which is homogeneously double G-linked [H-U] to $I_7 = I_6 :_R (y^{n - \beta}) = (x^{\alpha-\alpha_2},y^{\beta},z^{\gamma_1-\gamma_3},x^{\alpha_3 - \alpha_2}y^{\beta_3},x^{\alpha_1 - \alpha_2}y^{\beta_1},y^{\beta_2}z^{\gamma_2 - \gamma_3})$. Notice that $\mu(I_6) \leq 6$, thus $I_6$ is homogeneously glicci and thus $I$ is homogeneously glicci.
\\

\noindent\textbf{Subcase c):} If $\min\{\alpha_i\} = \alpha_3$, then $I_6 = (x^{\alpha-\alpha_3},y^{n},z^{\gamma-\gamma_3},y^{n - \beta + \beta_3},x^{\alpha_1 - \alpha_3}y^{n - \beta + \beta_1},\\x^{\alpha_2 - \alpha_3}y^{n - \beta}z^{\gamma_1 - \gamma_3},y^{n - \beta + \beta_2}z^{\gamma_2 - \gamma_3}) = (x^{\alpha-\alpha_3},y^{n},z^{\gamma-\gamma_3},x^{\alpha_1 - \alpha_3}y^{n - \beta + \beta_1},x^{\alpha_2 - \alpha_3}y^{n - \beta}z^{\gamma_1 - \gamma_3},\\y^{n - \beta + \beta_2}z^{\gamma_2 - \gamma_3})$. Notice that $\mu(I_6) \leq 6$, thus $I_6$ is homogeneously glicci and thus $I$ is homogeneously glicci.
\\

\noindent This concludes all possible cases, and thus we have proven that any $m$-primary monomial ideal $I \subset k[x,y,z]$ with $\mu(I) \leq 7$ is homogeneously glicci. We note that every second ideal in our construction is an $m$-primary monomial ideal.
\end{proof}

\noindent We may now finish the proof of our main result:

\begin{customthm}{3.1}{}

\noindent Let $R = k[x,y,z]$ be a polynomial ring over a field $k$. Let $I$ be an $m$-primary monomial ideal with $\mu(I) \leq 8$. Then $I$ is homogeneously glicci via a sequence of G-links with monomial ideals in each second step.
\end{customthm}

\begin{proof}
    
Let $I \subset k[x,y,z]$ be an $m$-primary monomial ideal with $\mu(I) \leq 8$. If $\mu(I) \leq 7$ then $I$ is homogeneously glicci by Lemma 3.2. Similarly, if $I$ is homogeneously licci, then $I$ is homogeneously glicci. If $I$ is not homogeneously licci, then $I$ is homogeneously linked to an ideal of the form $(x^\alpha,y^\beta,z^\gamma,x^{\alpha_1}y^{\beta_1},x^{\alpha_2}z^{\gamma_1},y^{\beta_2}z^{\gamma_2},x^{\alpha_3}y^{\beta_3}z^{\gamma_3},x^{\alpha_4}y^{\beta_4}z^{\gamma_4})$ by Theorem 2.1. Thus we may assume $I = (x^\alpha,y^\beta,z^\gamma,x^{\alpha_1}y^{\beta_1},x^{\alpha_2}z^{\gamma_1},y^{\beta_2}z^{\gamma_2},x^{\alpha_3}y^{\beta_3}z^{\gamma_3},x^{\alpha_4}y^{\beta_4}z^{\gamma_4})$. The links required turn out to be very similar to those needed in the seven generator case. First, start with $I_2 = ((x^{\alpha},y^{n},z^{\gamma}):_RI_1):_R(y^{n - \beta}).$ By Theorem 2.1, $I_2$ will be linked to $I_3 = ((x^{n},y^{n},z^{n}):_RI_2):_R(z^{n - \gamma}) = (x^{n},y^{n},z^{\gamma},x^{n - \alpha}y^{n - \beta + \beta_2}z^{\gamma_2}, x^{n - \alpha + \alpha_1}y^{n - \beta + \beta_1},x^{n - \alpha + \alpha_2}y^{n - \beta}z^{\gamma_1},x^{n - \alpha + \alpha_3}y^{n - \beta + \beta_3}z^{\gamma_3},\\x^{n - \alpha + \alpha_4}y^{n - \beta + \beta_4}z^{\gamma_4}).$ Next, by Corollary 2.3, $I_3$ is homogeneously G-linked to 
$I_4 = I_3 :_R (z^{\min(\gamma_i)}) = (x^{n},y^{n},z^{\gamma - \min(\gamma_i)},x^{n - \alpha}y^{n - \beta + \beta_2}z^{\gamma_2-\min(\gamma_i)}, x^{n - \alpha + \alpha_1}y^{n - \beta + \beta_1},x^{n - \alpha + \alpha_2}y^{n - \beta}z^{\gamma_1-\min(\gamma_i)},\\x^{n - \alpha + \alpha_3}y^{n - \beta + \beta_3}z^{\gamma_3 - \min(\gamma_i)},x^{n - \alpha + \alpha_4}y^{n - \beta + \beta_4}z^{\gamma_4-\min(\gamma_i)}).$ By Theorem 2.1, we have that $I_4$ is homogeneously G-linked to 
$I_5 = I_4 :_R (x^{n-\alpha}) = (x^{\alpha},y^{n},z^{\gamma - \min(\gamma_i)},y^{n - \beta + \beta_2}z^{\gamma_2-\min(\gamma_i)}, x^{\alpha_1}y^{n - \beta + \beta_1},\\x^{\alpha_2}y^{n - \beta}z^{\gamma_1-\min(\gamma_i)},x^{\alpha_3}y^{n - \beta + \beta_3}z^{\gamma_3 - \min(\gamma_i)},x^{\alpha_4}y^{n - \beta + \beta_4}z^{\gamma_4-\min(\gamma_i)}).$ By Theorem 2.1, we then have that $I_5$ will be homogeneously G-linked to
$I_6 = I_5 :_R (x^{\min(\alpha_i)}) = (x^{\alpha-\min(\alpha_i)},y^{n},z^{\gamma - \min(\gamma_i)},\\y^{n - \beta + \beta_2}z^{\gamma_2-\min(\gamma_i)}, x^{\alpha_1-\min(\alpha_i)}y^{n - \beta + \beta_1},x^{\alpha_2-\min(\alpha_i)}y^{n - \beta}z^{\gamma_1-\min(\gamma_i)},x^{\alpha_3-\min(\alpha_i)}y^{n - \beta + \beta_3}z^{\gamma_3 - \min(\gamma_i)},\\x^{\alpha_4-\min(\alpha_i)}y^{n - \beta + \beta_4}z^{\gamma_4-\min(\gamma_i)}).$ Finally, by Theorem 2.1, $I_6$ is homogeneously G-linked to
$I_7 = I_6 :_R (y^{n - \beta}) = (x^{\alpha-\min(\alpha_i)},y^{\beta},z^{\gamma - \min(\gamma_i)},y^{\beta_2}z^{\gamma_2-\min(\gamma_i)}, x^{\alpha_1-\min(\alpha_i)}y^{\beta_1},x^{\alpha_2-\min(\alpha_i)}z^{\gamma_1-\min(\gamma_i)},\\x^{\alpha_3-\min(\alpha_i)}y^{\beta_3}z^{\gamma_3 - \min(\gamma_i)},x^{\alpha_4-\min(\alpha_i)}y^{\beta_4}z^{\gamma_4-\min(\gamma_i)}).$\\

\noindent Thus if $\min(\alpha_i) \in \{\alpha_1,\alpha_2\}, \min(\beta_i) \in \min(\beta_1,\beta_2), \min(\gamma_i) \in \{\gamma_1,\gamma_2\}$, $\mu(I_7) \leq 7$ and thus $I$ is homogeneously glicci. We are left with the case where $\min\{\alpha_i\} \notin \{\alpha_1,\alpha_2\},  \min\{\beta_i\} \notin \{\beta_1,\beta_2\},\min\{\gamma_i\} \notin \{\gamma_1,\gamma_2\}$.\\

\noindent\textbf{Case 1:} Say $\min(\alpha_i) = \alpha_3, \min(\beta_i) = \beta_3,$ and $ \min(\gamma_i) = \gamma_3$. Then $I_7 = (x^{\alpha-min(\alpha_i)},y^{\beta_3},\\z^{\gamma - min(\gamma_i)},y^{\beta_2}z^{\gamma_2-min(\gamma_i)}, x^{\alpha_1-\min(\alpha_i)}y^{\beta_1},x^{\alpha_2-\min(\alpha_i)}z^{\gamma_1-min(\gamma_i)},x^{\alpha_4-min(\alpha_i)}y^{\beta_4}z^{\gamma_4-min(\gamma_i)})$. Thus since $\mu(I_7) \leq 7$, $I$ is homogeneously glicci. 
After re-labeling, this is the same as the $\alpha_4 = \min(\alpha_i)$ and $\gamma_4 = \min(\gamma_i)$ case.\\

\noindent\textbf{Case 2:} Say $\min(\alpha_i) = \alpha_3$ and $\min(\gamma_i) = \gamma_4$. For notational ease, re-write $I'_1 = I_7 =(x^{A},y^{B},z^{C},y^{B_2}z^{C_2}, x^{A_1}y^{B_1},x^{A_2}z^{C_1},y^{B_3}z^{C_3},x^{A_4}y^{B_4})$. Thus $I'_1$ is homogeneously G-linked to $I_2' = ((x^{A},y^{n},z^{C}):_RI_1'):_R(y^{n - b})$ by Theorem 2.1, which is homogeneously G-linked to $I_3' = ((x^{n},y^{n},z^{n}):_RI_2'):_R(z^{n - C}) = (x^{n},y^{n},z^{C},x^{n - A}y^{n - B + B_2}z^{C_2}, x^{n - A + A_1}y^{n - B + B_1},\\x^{n - A + A_2}y^{n - B}z^{C_1},x^{n - A + A_3}y^{n - B + B_3}z^{C_3}, x^{n - A + A_{4}}y^{n - B + B_{4}}z^{C_{4}}) = (x^{n},y^{n},z^{C},x^{n - A}y^{n - B + B_2}z^{C_2},\\ x^{n - A + A_1}y^{n - B + B_1},x^{n - A + A_2}y^{n - B}z^{C_1},x^{n - A}y^{n - B + B_3}z^{C_3}, x^{n - A + A_{4}}y^{n - B + B_{4}})$. Which, by Theorem 2.1, is homogeneously G-linked to $I'_4 = I'_3:_R(y^{n - B}) = (x^{n},y^{B},z^{C},x^{n - A}y^{B_2}z^{C_2}, x^{n - A + A_1}y^{B_1},\\x^{n - A + A_2}z^{C_1},x^{n - A}y^{B_3}z^{C_3}, x^{n - A + A_{4}}y^{B_{4}})$. By, Corollary 2.3, $I'_4$ is homogeneously double G-linked to $I'_5 = I'_4:_R(y^{\min(B_i)}) = (x^{n},y^{B - \min(B_i)},z^{C},x^{n - A}y^{B_2 - \min(B_i)}z^{C_2}, x^{n - A + A_1}y^{B_1-\min(B_i)},\\x^{n - A + A_2}z^{C_1},x^{n - A}y^{B_3 - \min(B_i)}z^{C_3}, x^{n - A + A_{4}}y^{B_{4}-\min(B_i)}).$ Now, we notice that if $\min(B_i) = B_4$, we are done. For our final case, if $\min(B_i) = B_3$, then $ I'_5 = I'_4:_R(y^{B_3}) = (x^{n},y^{B - B_3},z^{C},\\x^{n - A}y^{B_2 - B_3}z^{C_2},x^{n - A + A_1}y^{B_1-B_3},x^{n - A + A_2}z^{C_1},x^{n - A}z^{C_3}, x^{n - A + A_{4}}y^{B_{4}-B_3}).$ $I_5'$ is homogeneously G-linked to $ I'_6 = I'_5:_R(x^{n - A}) = (x^A,y^{B - B_3},z^{C},y^{B_2 - B_3}z^{C_2}, x^{A_1}y^{B_1-B_3},x^{A_2}z^{C_1},z^{C_3}, x^{A_{4}}y^{B_{4}-B_3})$ and clearly $\mu(I'_6) \leq 7.$ After re-labeling, this is the same as the $\alpha_4 = \min(\alpha_i)$ and $\gamma_3 = \min(\gamma_i)$ case.
Thus $\mu(I) = 8$ implies $I$ homogeneously glicci. 
\end{proof}

\noindent Notice that throughout this paper, after each homogeneous double G-link we are left with another monomial ideal. This is a rare property when studying homogeneously glicci ideals, as the only Gorenstein monomial ideals are pure-power ideals (which are complete intersections). The difficulty in expanding our result to $m$-primary monomial ideals with larger minimal generating sets lies in $\mu(K_{xy}),\mu(K_{xz}),$ and $\mu(K_{yz})$. When $\mu(I) \geq 9$, it is possible for $\mu(K_{xy}),\mu(K_{xz}),$ and $\mu(K_{yz})$ to all be $\geq 2$. Any $m$-primary monomial ideal that is homogeneously linked to such an ideal will be nonlicci and Corollary 2.3 cannot be applied for further reduction. It seems likely that linking such ideals to a complete intersection would require the use of double G-links that link to non-monomial ideals. Despite this limitation, our methods may still be used to obtain a large class of homogeneously glicci (but nonlicci) $m$-primary monomial ideals in a more general context.

\section{General Construction}
We now provide a recipe for constructing homogeneously glicci $m$-primary monomial ideals in $R = k[x_1,...,x_n]$ with no restrictions on the sizes of their minimal generating sets. Let $I_0$ be any $m$-primary monomial ideal in $R$ such that $I_0^{\#}$ is generated by monomials in $k[x_1,x_2,x_3]$ and $\mu(I_0) \leq n + 5$. Write $I_0 = (x_1^{\alpha^{(0)}_{1}},...,x_n^{\alpha^{(0)}_{n}} )+I_0^{\#}$. For $i > 0$, choose integers $\alpha^{(i)}_{1},...,\alpha^{(i)}_{n}$ such that for all $j$, $\alpha_j^{(i)} > \alpha_j^{(i-1)}$. Take $\beta_j^{(i)} < \alpha_j^{(i)} - \alpha_j^{(i-1)}$ and define $I_i = x_1^{\beta^{(i)}_1}\cdots x_n^{\beta_n^{(i)}}(I_{i-1}) + (x_1^{\alpha^{(i)}_1},...,x_n^{\alpha^{(i)}_n})$. Continuing this process yields a sequence $I_0,I_1,...$ of $m$-primary glicci monomial ideals. When $\rm{ht}(I^{\#}_0) \geq 2$, they are not licci.

\begin{proof}
    Say $J = I_k$ for some $k \in \mathbb{N}$. We know that $I_k$ is homogeneously linked to $I_{k-1}$ by Theorem 2.1. Continuing this process shows that $J$ is homogeneously linked to $I_0$, which is homogeneously glicci by Theorem 3.1. Thus $J$ is homogeneously glicci. If ${\rm ht}(I_0^{\#})\geq 2$, then $J$ is not licci by Theorem 2.1.    
\end{proof}

\printbibliography

\end{document}